\documentclass[a4paper,12pt]{amsart}

\usepackage{amssymb,amsbsy,amsmath,amsfonts,amssymb,amscd}
\usepackage{latexsym}
\usepackage{graphics}
\usepackage{url,hyperref}

\usepackage{euscript,color}

\def\LL{\mathcal{L}}

\hypersetup{
    bookmarks=true,         
    unicode=false,          
    pdftoolbar=true,        
    pdfmenubar=true,        
    pdffitwindow=false,     
    pdfstartview={FitH},    
    pdftitle={My title},    
    pdfauthor={Author},     
    pdfsubject={Subject},   
    pdfcreator={Creator},   
    pdfproducer={Producer}, 
    pdfkeywords={keyword1} {key2} {key3}, 
    pdfnewwindow=true,      
    colorlinks=true,       
    linkcolor=black,          
    citecolor=green,        
    filecolor=magenta,      
    urlcolor=blue         
}

\newtheorem{theorem}{Theorem}[section]

\newtheorem{lemma}[theorem]{Lemma}
\newtheorem{prop}[theorem]{Proposition}

\newtheorem{thm}{Theorem}[section]

\begin{document}

\title[Killing vector fields of constant length]{Killing vector fields of constant length on compact hypersurfaces }
\author{ Antonio J. Di Scala}

\date{\today}

\maketitle

{\bf Abstract:} We show that if a compact hypersurface $M \subset \mathbb{R}^{n+1}$, $n \geq3$, admits a non zero Killing vector field $X$ of constant length then $n$ is even and $M$ is diffeomorphic to the unit hypersphere of $\mathbb{R}^{n+1}$.
Actually, we show that $M$ is a complex ellipsoid in $\mathbb{C}^{N} = \mathbb{R}^{n+1}$.
As an application we give a simpler proof of a recent theorem due to S. Deshmukh \cite{De12}.

\section{Introduction}
In a recent paper \cite{De12}, S. Deshmukh motivated by the Hopf vector fields of odd dimensional hyperspheres
asked the following question:\\

\emph{ ``Does an orientable compact hypersurface in a Euclidean space that admits a unit Killing vector field $\xi$ together with the condition that $g(A\xi,\xi)$ is a constant necessarily isometric to a sphere?[sic]"}\\

Deshmukh answer positively but his proof is not correct. In \cite[page 972]{De12} after equation (7) he claims that the shape operator $A$ of a hypersurface $M$ satisfies \[ d \psi_t \circ A = A \circ d \psi_t  \]
where $\{\psi_t \}$ is the flow of the unit Killing vector field $\xi$.
This claim is false since the shape operator $A$ contains extrinsic information about $M$.
For example, if  $M$ is a plane curve and $t$ its arc length parameter then $\xi = \partial_t$ is a unit Killing vector field.
But then the above equation implies that the curvature (i.e. the inverse of the radius of the osculating circle) of the plane curve $M$ is constant which is not true for all plane curves.

It is known that compact hypersurfaces are orientable (see for example \cite{Sa69}) so such hypothesis can be omitted from the above question.\\

In this article we show that besides Deskmush's wrong proof the answer to his question is indeed positive.
Moreover, we show that the condition $g(A\xi,\xi)$ to be constant, where $g$ is the induced Riemannian metric on $M$,  is not necessary to characterize $M$ topologically and geometrically.\\

 Here is our main result.\\

\begin{thm} \label{introduction} A compact hypersurface $M \subset \mathbb{R}^{n+1}$  that admits a unit Killing vector field $\xi$ is diffeomorphic to an odd dimensional hypersphere. Actually, $n+1 = 2N$ and  $M$ is (up to isometry $\mathbb{R}^{n+1} \cong \mathbb{C}^{N}$)  the complex ellipsoid defined by the equation \[ |\omega_1 \, z_1|^2 + |\omega_2  \, z_2|^2 + \cdots + |\omega_N \, z_N|^2 = 1 \, \]
where $(z_1,\cdots,z_N)$ are complex coordinates of $\mathbb{C}^N $ and $\omega_1,\cdots,\omega_N$ are non zero real numbers. Moreover,  if  $g(A \xi,\xi)$ is constant then \[ \omega_1^2 = \cdots = \omega_N^2 \] and $M$ is a hypersphere of constant curvature.
\end{thm}

For more about Killing vector fields of constant length on Riemannian manifolds see \cite{BN08}.

\section{Preliminaries}

A Killing vector field $X$ of the flat Euclidean space $\mathbb{R}^{n+1}$ is given by a pair $X = (R ,\mathbf{v})$ where $R$ is a skew-symmetric matrix and $\mathbf{v} \in \mathbb{R}^{n+1}$ is a vector.
At a point $\mathbf{p} \in \mathbb{R}^{n+1}$, $X(\mathbf{p})$ is given by the formula \[ X(\mathbf{p}) = R \, \mathbf{p} + \mathbf{v} \, .\]


The block diagonal form for skew-symmetric matrices can be used to introduce complex coordinates $(z_1,\cdots,z_N)$ on the image of $R$. So, we have that $\mathbb{R}^{n+1} = \mathbb{C}^{N} \oplus  \mathbb{R}^{k}$ and
the Killing vector field $X$ is given by \begin{equation}\label{killing}X(z_1, \cdots, z_N, \mathbf{r})  = (i \omega_1 z_1 , \cdots, i \omega_N  z_N, \mathbf{v}_k) \end{equation}
where $\mathbf{v}_k$ is the projection of $\mathbf{v}$ to $\mathbb{R}^k$ and $\omega_1,\cdots,\omega_N$ are non zero real numbers.

Let $H_X^{\LL} \subset \mathbb{R}^{n+1} $, $\LL > 0$, be the set of points where the Killing vector field $X$ has constant length ${\LL}$. The following is a straightforward consequence of equation (\ref{killing}).

\begin{prop}\label{conexas} The set $H_X^{\LL} \subset \mathbb{R}^{n+1}$ is a smooth connected submanifold and the Killing vector field $X$ is always tangent to $H_X^{\LL}$. Moreover, if $H_X^{\LL} \neq \emptyset$ then one of the following holds:
\begin{itemize}
\item[(i)]  $\LL = \| \mathbf{v}_{ker}\|$ and $H_X^{\LL} = \mathbb{R}^{k}$ ,
\item[(ii)]  $\LL > \| \mathbf{v}_{ker}\|$ and  $H_X^{\LL} = \mathrm{E} \times \mathbb{R}^k \subset \mathbb{C}^N \times \mathbb{R}^k$ where  $\mathrm{E} \subset \mathbb{C}^N$ is the complex ellipsoid defined by \[ |\omega_1 z_1|^2 + \cdots + |\omega_N z_N|^2 = \LL - \| \mathbf{v}_{ker}\| \, ,\]
    where $\omega_1,\cdots,\omega_N$ are non zero real numbers.
\end{itemize}
\end{prop}

Here is another important consequence of equation (\ref{killing}).

\begin{lemma}\label{type} The rank of the shape operator $S_p$ of the hypersurface $H_X^{\LL}$ is constant i.e. it does not depends upon the point $p \in H_X^{\LL}$. Moreover, in case {\rm (ii)} of the previous proposition the rank of $S$ is $2N-1$ at all points of  $H_X^{\LL}$.
\end{lemma}
\it Proof. \rm Indeed, the vector field $\mathbf{n} = (\omega_1^2 \, z_1 , \cdots, \omega_N^{2} z_N)$ is perpendicular to the complex ellipsoid $\mathrm{E}$ at $p=(z_1,\dots,z_N)$. Then for $\mathbf{t} \in \mathrm{T}_p \mathrm{E}$ we have \[ S_p(\mathbf{t}) . \mathbf{t} = -D_\mathbf{t} \mathbf{n}.\mathbf{t} = -\sum_{j=1}^N \omega_j^2 |t_j|^2 \, .\]
Since $\omega_1,\cdots,\omega_N$ are non zero real numbers we get that the rank of the shape operator $S_p$ is $2N-1$. $\Box$

\section{Main theorem}

\begin{thm}\label{rigid} Let $(M^n,g)$ be a connected compact Riemannian manifold of dimension $n \geq 3$.
Let $\xi$ be a non zero Killing vector field of $M$ with constant length $\LL = \|\xi\| $.\\
If $(M,g)$ can be isometrically embedded into the flat Euclidean space $\mathbb{R}^{n+1}$ then $n$ is even and
\begin{itemize}
 \item[$(i)$] the Killing vector field $\xi$ extends to a Killing vector field $X = (R,\mathbf{v})$ of $\mathbb{R}^{n+1}$,
 \item[$(ii)$] $M = H^{\LL}_{X}$.
\end{itemize}
\end{thm}

\it Proof. \rm  We will identify $M$ with its image under the isometric embedding into $\mathbb{R}^{n+1}$.
The idea is that for $n \geq 3$ hypersurfaces of $\mathbb{R}^{n+1}$ whose shape operator has rank $\geq 3$ are rigid due to the theorem of Beez-Killing \cite[page 46, Corollary 6.5]{KN69}.

Let $A$ be the shape opertator of $M$. It is well known that $M$ has a point $q$ where $A_q$ is positive (or negative) definite \cite[page 255]{GHL04}. Since the rank of $A_q$ is $n$ the theorem of Beez-Killing implies that the Killing vector field $\xi$ extends to a Killing vector field $X$ of $\mathbb{R}^{n+1}$ near $q$. Then near $q$ the hypersurface $M$ is contained in  $H^{\LL}_{X}$. If $\mathrm{dim}(H^{\LL}_{X}) = n+1$ then
$X$ is a parallel vector field i.e. $X = \mathbf{v}$ a constant vector. Then $\xi$ is the restriction of a constant vector and so $A_q(\xi) = 0$ which contradicts that $A_q$ is definite at $q$. So $\mathrm{dim}(H^{\LL}_{X}) = n$ and $M = H^{\LL}_{X}$ near $q$. Notice that $H^{\LL}_{X}$ is as in {\rm (ii)} of Proposition \ref{conexas}.
That is to say, the shape operator of  $H^{\LL}_{X}$ has rank $n= 2N - 1$ at $q$ which show that $n$ is even.\\

Let $U = M \bigcap H^{\LL}_{X}$ be the intersection of both hypersurfaces and let $\stackrel{\circ}{U}$ be the interior of $U$ in $M$. By the previous observation $\stackrel{\circ}{U}$ contains a neighborhood of $q$ in $M$ an so $\stackrel{\circ}{U}$ is non empty. We claim that $\stackrel{\circ}{U}$ is closed in $M$. Indeed, let $p_j \in \, \stackrel{\circ}{U}$ be a sequence converging to $p_0 \in M$.  The shape operator $A_{p_0}$ is the limit of the shape operator $S_{p_j}$ of $H^{\LL}_{X}$. Then $A_{p_0} = S_{p_0}$ and Lemma \ref{type} implies that the rank of $A_{p_0}$ is $n \geq 3$. Again by the theorem of  Beez-Killing we have that near $p_0$ the vector field $\xi$ extends to a Killing vector field $Y$ of $\mathbb{R}^{n+1}$. Notice that $Y$ restricted to $\stackrel{\circ}{U}$ agree with $X$ near the points $p_j$ for $j>>0$. This implies $X=Y$ on $\mathbb{R}^{n+1}$  and then  $p_0 \in \, \, \stackrel{\circ}{U}$. So $\stackrel{\circ}{U} \, \, = M$ hence $M \subset H^{\LL}_{X}$. This shows $(i)$ and $(ii)$ since by Proposition \ref{conexas} the compact hypersurface $H^{\LL}_{X}$ is connected. $\Box$

\subsection{Proof of Theorem \ref{introduction}}
It is clear that the first part of Theorem \ref{introduction} follows from Theorem \ref{rigid}.
By using the shape operator as in the proof of Lemma \ref{type} the condition $g(A \xi,\xi)$ is constant
gives:  \[\omega_1^4 \,|z_1|^2 + \cdots + \omega_N^4 |z_N|^2 = cte \, ,\]
on the complex ellipsoid $|\omega_1 \, z_1|^2 + \cdots + |\omega_N \, z_N|^2 = 1$.
Then \[ \omega_1^4 \,|z_1|^2 + \cdots + \omega_N^4 |z_N|^2 = cte \left( |\omega_1 \, z_1|^2 + \cdots + |\omega_N \, z_N|^2 \right) \]
for all $(z_1,\cdots,z_N) \in \mathbb{C}^N$. So $\omega_j^2 = cte$ for $j=1,\cdots,N$ and $M$ is the hypersphere of $\mathbb{C}^N$defined by the equation
\[ |z_1|^2 + \cdots + | z_N|^2 = \frac{1}{cte} \, . \]
\hspace{12cm} $\Box$

\section{Further remarks}
Surfaces in $\mathbb{R}^3$ are not rigid so it is natural to ask about existence of compact surfaces with a unit Killing vector field in $\mathbb{R}^3$. The answer is that such surface does not exists. Indeed, it is not difficult to see that a unit Killing vector field of a surface is a parallel vector field. So the metric is flat. But a flat compact surface cannot be a surface in $\mathbb{R}^3$ due to the classical theorem about existence of elliptic points \cite[page 255]{GHL04}.

\vspace{1cm}

\noindent Antonio J. Di Scala, \\ Dipartimento di Scienze Matematiche  ``G. L. Lagrange", \\
Politecnico di Torino, Corso Duca degli Abruzzi 24, \\
10129 Torino, Italy. \\
    email:     antonio.discala@polito.it \\

\end{document}